\documentclass[preprint]{elsarticle}

\usepackage{cmap} 
           \usepackage[T2A]{fontenc}   
           \usepackage[utf8]{inputenc}
           \usepackage[english]{babel}
\usepackage{amsmath}
\usepackage{amsfonts}
\usepackage{amssymb}
\usepackage{amsthm}
\usepackage[pdftex,unicode]{hyperref}
\usepackage{indentfirst} 
\usepackage{amscd}
\usepackage[all,cmtip]{xy}
\usepackage{mathtools}

\def\om{\omega}

\def\al{\alpha}

\def\la{\lambda}

\def\nfs/{NFS}
\def\cdp/{CDP}
\def\cdpz/{CDP${}_0$}

\def\cl#1{\overline{#1}}

\def\es{\varnothing}
\def\Tau{{\mathcal T}}

\def\set#1{\bbset#1\eeset}
\def\bbset#1:#2\eeset{\{#1\,:\,#2\}}

\def\bbsett#1:#2\eesett{\{#1\,:\,\text{#2}\}}

\def\ibbset#1:#2\ieeset{(#1)_{#2}}

\def\gs{{\mathfrak{S}}}

\def\tp{\Tau}

\def\cP{{\mathcal P}}
\def\cB{{\mathcal B}}

\def\cB{{\mathcal B}}

\newcommand\restrA[2]{{
  \left.\kern-\nulldelimiterspace 
  #1 
  \vphantom{\big|} 
  \right|_{#2} 
  }}

\def\pwr#1_#2{#1^{[#2]}}

\def\D{\Delta}
\def\term#1{{\it #1}}

\def\dddm#1(#2){N_{#1}(#2)}
\def\dddb#1(#2){B_{#1}(#2)}

\def\et(#1){ (#1)}

\def\bitem#1,#2.{ $#2\nrightarrow #1$:\ }

\newtheorem{proposition}{Proposition}
\newtheorem{theorem}{Theorem}
\newtheorem*{theorem*}{Theorem}

\newtheorem*{lemma*}{Lemma}
\newtheorem{cor}{Corollary}

\theoremstyle{definition}

\theoremstyle{remark}

\newtheorem*{note*}{Remark}

\def\oo#1/{$O_{#1}$}

\def\gd/{$G_\delta$}

\def\rarr{\Rightarrow}

\def\cnt{\Lambda}

\def\pch{\pi\chi}
\def\gs{\mathfrak{s}}

\begin{document}

\begin{frontmatter}

\title{Metrizability of CHART groups}
\author{Evgenii Reznichenko} 

\ead{erezn@inbox.ru}

\address{Department of General Topology and Geometry, Mechanics and  Mathematics Faculty, M.~V.~Lomonosov Moscow State University, Leninskie Gory 1, Moscow, 199991 Russia}

\begin{abstract}
For compact Hausdorff admissible right topological (CHART) group $G$, we prove $w(G)=\pi\chi(G)$.
This equality is well known for compact topological groups.
This implies the criteria for the metrizability of CHART groups: if $G$ is first-countable (2013, Moors, Namioka) or $G$ is Fr\'echet (2013, Glasner, Megrelishvili), or $G$ has countable $\pi$-character (2022, Reznichenko) then $G$ is metrizable. Under the continuum hypothesis (CH) assumption, a sequentially compact CHART group is metrizable.
Namioka's theorem that metrizable CHART groups are topological groups extends to CHART groups with small weight.
\end{abstract}
\begin{keyword}
compact right topological groups
\sep
admissible  groups
\sep
CHART groups
\sep
metrizable spaces
\sep
$\pi$-character
\end{keyword}
\end{frontmatter}

\section{Introduction}

A group $G$ with a topology is called \term{right topological} if all right shifts $\rho_h: G\to G, g\mapsto gh$ are continuous. The set of  $g\in G$ for which the left shift $\la_g: G \to G, h\mapsto gh$ is continuous is called the \term{topological center} and is denoted as $\cnt(G)$.
A group $G$ with topology is called \term{semitopological} if all right and left shifts are continuous, that is, if $G$ is right topological and $\cnt(G)=G$.
A right topological group $G$ is called \term{admissible} if $\cnt(G)$ is a dense subset of $G$.
We write ``CHART'' for ``compact Hausdorff admissible right topological''.
%
A \term{paratopological group} $G$ is a group $G$ with a topology such that the product map of $G \times G$ into $G$ is jointly continuous.

The study of groups with topology in which not all operations are continuous and conditions implying the continuity of operations began with the 1936 paper \cite{Montgomery1936} of Montgomery, who, among other things, proved that a Polish (i.e., separable metrizable by a complete metric) semitopological group is a topological group. Interest in such groups was renewed in relation to topological dynamics. The autohomemorphism group of a locally compact space in the compact open topology is a paratopological group \cite{Arens1946}. In the same paper Arens obtained conditions under which an autohomemorphism group is a topological group. In 1957 Ellis proved that a locally compact paratopological group is a topological group \cite{Ellis1957}. Shortly afterwards,  he strengthened this theorem to semitopological locally compact groups and proved the celebrated Ellis theorem \cite{Ellis1957-2}: any semitopological locally compact group is a topological group. The enveloping semigroup of a dynamical system was introduced by Ellis in 1960 \cite{Ellis1960}. It has become a fundamental tool in the abstract theory of topological dynamical systems. The study of dynamical systems for which the enveloping semigroup is a CHART group plays a large role in the abstract theory of topological dynamical systems.

Let $G$ be a CHART group. Any of the following conditions implies that $G$ is a topological group.
\begin{itemize}
\item[$(C_1)$]
$G$ is metrizable (Theorem 2.1 \cite{namioka1972}).
\item[$(C_2)$]
$G$ is first-countable (Remark after Proposition 2.7 \cite{MoorsNamioka2013}).
\item[$(C_3)$]
$G$ is Fr\'echet (Corollary 8.8 \cite{GlasnerMegrelishvili2013}).
\item[$(G_4)$]
$G$ has countable $\pi$-character (for example, $G$ is a compact space with countable tightness) (Corollary 2 (3) \cite{rezn2022-1}).
\end{itemize}

A CHART group $G$ is \term{tame} if for every $g \in G$, the mapping $x \mapsto g \cdot x$ is fragmented.
A compact semitopological group is tame.

Note even more conditions that imply that $G$ is a topological group:
$(G_5)$ the multiplication of $G$ is separately continuous (Ellis theorem \cite{Ellis1957-2});
$(G_6)$ $G$ is tame \cite[Theorem 21]{GlasnerMegrelishvili2013};
$(G_7)$ the multiplication of $G$ is continuous at $(e,e)$; (follows from \cite{ruppert1975} and \cite[Theorem 5]{Milnes1979});
$(G_8)$ the multiplication of $G$ is feebly continuous \cite[Proposition 3.2]{Moors2016} (see also \cite[Corollary 2 (2)]{rezn2022-1});
$(G_9)$ the inversion $g\mapsto g^{-1}$ is continuous at $e$ \cite[Theorem 5]{Milnes1979};
$(G_{10})$ $\cnt(G)$ is a topological group (or merely contains a dense topological group) \cite[Theorem 5]{Milnes1979};
the right translations $g\mapsto gh$ form an equicontinuous family of maps from $G$ onto $G$ \cite[Theorem 5]{Milnes1979}.

Note that the group $G$ in $(C_1)$--$(C_4)$ is metrizable because compact first-countable, Fr\'echet topological groups and groups of countable $\pi$-character are metrizable (Corollary 4.2 .2 and Corollary 5.7.26 of \cite{at2009}). Recall that for compact Hausdorff spaces \cite{Hodel1984handbook}:
\begin{align*}
 \text{metrizable}\rarr
 \text{first-countable}\rarr
 \text{Fr\'echet}\rarr
\\
 \text{countable tightness}\rarr
 \text{countable $\pi$-character.}
\end{align*}

In this note, we prove that $w(G)=\pch(G)$ for CHART group $G$ (Theorem \ref{t:pichi:2}). Whence it follows that $G$ is metrizable if $G$ has countable $\pi$-character. Note that this fact reduces $(C_2)$--$(C_4)$ to $(C_1)$.

Under the continuum hypothesis (CH) assumption, a sequentially compact CHART group is metrizable (Corollary \ref{c:seqcomp:1}).
Namioka's theorem that metrizable CHART groups are topological groups extends to CHART groups with small weight (Theorem \ref{t:ma:2}).

\section{Definitions and notation}

We will denote by $G$ a group and by $e\in G$ the identity of the group.

Let $X$ be a space, $x\in X$, $\cP$ be a family of open subsets of $X$.
A family $\cP$ is called a \term{base in $x$} if for any neighborhood $U$ of $x$ there exists $V\in \cP$ so that $x\in V\subset U$.
A family $\cP$ is called a \term{$\pi$-base in $x$} if for any neighborhood $U$ of $x$ there exists $V\in \cP$ so that $V\subset U$.
A family $\cP$ is called a \term{base of $X$} if $\cP$ is a base at every point of $X$.
A family $\cP$ is called a \term{$\pi$-base of $X$} if $\cP$ is a $\pi$-base at every point of $X$.
Denote the diagonal $\D_X=\set{(x,x):x\in X}$ in $X^2$.

Recall the necessary definitions of cardinal functions from \cite{Hodel1984handbook}.
\begin{itemize}
\item[]
\term{weight}
\[
w(X)=\min\set{|\cB|: \cB\text{ a base for }X};
\]
\item[]
\term{diagonal degree}
\begin{align*}
\D(X)&=\min\{ |\cP|: \cP\text{ a family of open neighborhoods}
\\
&\text{of the diagonal }\D_X\text{ and }\bigcap \cP=\D_X\};
\end{align*}
\item[]
\term{character}
\begin{align*}
\chi(x,X)&=\min\{ |\cP|: \cP\text{ a base for }x\},
\\
\chi(X)&=\sup\{ \chi(x,X): x\in X\};
\end{align*}
\item[]
\term{$\pi$-character}
\begin{align*}
\pch(x,X)&=\min\{ |\cP|: \cP\text{ a $\pi$-base for }x\},
\\
\pch(X)&=\sup\{ \pch(x,X): x\in X\};
\end{align*}
\item[]
\term{tightness}
\begin{align*}
t(x,X)=\min\{& \tau: \text{ for all }A\subset X\text{ with }x\in\cl A
\\
&\text{there is }M \subset A\text{ with }|M|\leq\tau\text{ and }x\in\cl M\},
\\
t(X)&=\sup\{ t(x,X): x\in X\}.
\end{align*}
\end{itemize}

A space $X$ is called \term{sequentially compact} provided that every sequence in $X$
has a convergent subsequence.

A topological space $X$ is said to satisfy the \term{countable chain condition}, or to be \term{ccc}, if the partially ordered set of non-empty open subsets of $X$ satisfies the countable chain condition, i.e. every pairwise disjoint collection of non-empty open subsets of $X$ is countable.

In what follows, it is assumed that the spaces are Hausdorff.

\section{Diagonal degree of groups}\label{sec-pichi}

\begin{theorem}\label{t:pichi:1}
Let $G$ be a right topological group and $\cnt(G)^{-1}$ dense in $G$.
Then $\D(G)\leq \pch(G)$.
\end{theorem}
\begin{proof}
Let $\tau=\pch(G)$ and $\set{U_\al:\al<\tau}$ be a $\pi$-base for $e$.
We set $C=\cnt(G)$ and
\[
W_\al=\bigcup_{g\in C} gU_\al \times gU_\al
\]
for $\al<\tau$. Let us show that $\bigcap_{\al<\tau}W_\al=\D_X$.

Let us show that $\D_X\subset \bigcap_{\al<\tau}W_\al$.
Let $x\in G$ and $\al<\tau$. Since $C^{-1}$ is dense in $G$, then $g^{-1}\in U_\al x^{-1}$ for some $g\in C$. Then
\[
(x,x)\in gU_\al \times gU_\al \subset W_\al.
\]
Let us show that $\bigcap_{\al<\tau}W_\al \subset \D_X$. Assume the opposite, that is, there is
\[
(x,y)\in \bigcap_{\al<\tau}W_\al \setminus \D_X.
\]
Since $x\neq y$ the group $G$ is Hausdorff and right topological, there exists a neighborhood $U$ of the identity for which $Ux^{-1}\cap Uy^{-1}=\es$. Since $\set{U_\al:\al<\tau}$ is a $\pi$-base for $e$, then $U_\al\subset U$ for some $\al<\tau$. There is $g\in C$ so $(x,y)\in gU_\al \times gU_\al$. Then $x,y\in gU_\al$ and
$g\in U_\al x^{-1}\cap U_\al y^{-1}$. Hence $g\in Ux^{-1}\cap Uy^{-1}$, a contradiction.
\end{proof}

Theorem \ref{t:pichi:1} is a generalization of \cite[Corollary 2.5]{arh-rezn2005}.

\begin{proposition}\label{p:pichi:1}
If $G$ is a CHART group, then $\cnt(G)=\cnt(G)^{-1}$ is a subgroup.
\end{proposition}
\begin{proof}
Clearly, $\cnt(G)$ is a subsemigroup.
Let $g\in\cnt(G)$.
Since $\la_g$ is a continuous bijection of a compact space, then $\la_g$ is a homeomorphism and the mapping $\la_g^{-1}=\la_{g^{-1}}$ is continuous.
\end{proof}

\begin{theorem}\label{t:pichi:2}
Let $G$ be a CHART group.
Then $w(G)= \pch(G)$.
\end{theorem}
\begin{proof}
Always $w(G)\leq \pch(G)$. Proposition \ref{p:pichi:1} and Theorem \ref{t:pichi:1} imply $\D(G)\leq \pch(G)$.
For compact Hausdorff spaces $w(G)=\D(G)$ \cite[Corollary 7.6]{Hodel1984handbook}.
Hence $w(G)= \pch(G)$.
\end{proof}

\begin{cor}\label{c:pichi:1}
Let $G$ be a CHART group.
Then the following conditions are equivalent.
\begin{itemize}
\item[$(1)$]
$G$ is metrizable;
\item[$(2)$]
$G$ is first-countable;
\item[$(3)$]
$G$ is Fr\'echet;
\item[$(4)$]
$G$ has countable tightness;
\item[$(5)$]
$G$ has countable $\pi$-character.
\end{itemize}
\end{cor}
\begin{proof}
Always $(1)\rarr(2)\rarr(3)\rarr(4)$.
Since for compact spaces $\pch(G)\leq t(G)$ \cite[Theorem 7.13]{Hodel1984handbook}, $(4)\rarr(5)$ is true. Theorem \ref{t:pichi:2} implies $(5)\rarr(1)$.
\end{proof}

\section{Sequentially compact CHART groups}\label{sec-seqcomp}

Denote $I=[0,1]$,
\[
\gs=\min \set{\tau: I^\tau \text{ is not sequentially compact}}.
\]
The cardinal $\gs$ is called \term{splitting number}, $\om < \gs \leq 2^\om$ \cite{vanDouwen1984handbook,Vaughan1984handbook}.

\begin{proposition}\label{p:seqcomp:1}
Let $X$ be a compact sequentially compact space.
Then 
\[
\pch(x,X)<\gs
\]
for some $x\in X$.
\end{proposition}
\begin{proof}
Let us assume the opposite. Then $\pch(x,X)\geq\gs$ for all $x\in X$. Then \cite[Theorem 1]{Shapirovskii1980} implies that $X$ maps continuously onto $I^\gs$. Since sequential compactness is preserved by continuous mappings, $I^\gs$ is a sequentially compact space. Contradiction.
\end{proof}

If $G$ is a compact Hausdorff space and $w(G)=\tau<\gs$ then $G$ can be embedded in $I^\tau$ and $G$ is sequentially compact.
In a homogeneous space $G$, if at some point the $\pi$-character is equal to $\tau$, then the $\pi$-character of the whole space is equal to $\tau$. Therefore, the Proposition \ref{p:seqcomp:1} and Theorem \ref{t:pichi:2} imply the following proposition.

\begin{theorem}\label{t:seqcomp:1}
Let $G$ be a CHART group.
Then $G$ is sequentially compact if and only if $w(G) < \gs$.
\end{theorem}

Assuming the conitinuum hypothesis (CH), $\om < \gs \leq 2^\om=\om_1$, that is, $\gs=\om_1$. Theorem \ref{t:seqcomp:1} implies the following proposition.

\begin{cor}\label{c:seqcomp:1}
{\rm(CH)}
Let $G$ be a sequentially compact CHART group.
Then $G$ is metrizable.
\end{cor}

\section{Martin's axiom and continuity of operations in CHART groups}\label{sec-ma}

Recall the topological characterization of the statement $MA(\tau)$:
\begin{itemize}
\item[$MA(\tau)$]
if $X$ is a compact Hausdorff topological space that satisfies the ccc then $X$ is not the union of $\tau$ or fewer nowhere dense subsets.
\end{itemize}
\term{Martin's axiom} (MA): For every $\tau< 2^\om$, $MA(\tau)$ holds.

A topological space $(X,\tp)$ is called \term{$\D_s$-nonmeager} space \cite{rezn2022gbd,rezn2022-1} if for any mapping $\Omega: X\to\tp$ such that $ x\in \Omega(x)$ for $x\in\ X$, there exists a nonempty $W\in\tp$, such that
\begin{equation}\label{eq:ma:1}
W\subset\cl{\set{x\in W: W\subset \Omega(x)}}.
\end{equation}

\begin{proposition}\label{p:ma:1}
$MA(\tau)$.
Let $X$ be a ccc compact space and $w(G)\leq\tau$. Then $X$ is $\D_s$-nonmeager.
\end{proposition}
\begin{proof}
Let $\tp$ be the topology of $X$,
$\set{U_\al:\al<\tau}\subset \tp$ is the base of $X$ and $\Omega: X\to\tp$ is a mapping such that $x\in \Omega(x) $ for $x\in\ X$.
For $\al<\tau$, put $M_\al=\set{x\in U_\al: U_\al\subset \Omega(x)}$. Then $G=\bigcup_{\al<\tau}M_\al$. $MA(\tau)$ implies that $\cl{M_\al}$ has nonempty interior for some $\al<\tau$. Take a nonempty open $W\subset U_\al$ such that $M_\al \cap W$ dense in $W$. Then (\ref{eq:ma:1}) holds.
\end{proof}

Proposition \ref{p:pichi:1} and \cite[Theorem 13 and Theorem 17(2)]{rezn2022-1} imply the following assertion.

\begin{theorem}\label{t:ma:1}
If $G$ is a $\D_s$-nonmeager CHART group then $G$ is topological group.
\end{theorem}

\begin{theorem}\label{t:ma:2}
$MA(\tau)$.
Let $G$ be a CHART group.
If $w(G)\leq\tau$ then $G$ is a topological group.
\end{theorem}
\begin{proof}
CHART groups have a right-invariant Haar measure \cite{MilnesPym1992,Moors2015}. Hence $G$ is a ccc space.
It follows from Proposition \ref{p:ma:1} that $G$ is $\D_s$-nonmeager.
Theorem \ref{t:ma:1} implies that $G$ is a topological group.
\end{proof}

\begin{cor}[{Corollary 3 \cite{rezn2022-1}}]\label{c:ma:1}
{\rm (MA)}
Let $G$ be a CHART group.
If $w(G)< 2^\om$ then $G$ is a topological group.
\end{cor}

Since $MA(\om)$ is true in ZFC, Theorem \ref{t:ma:2}
the following assertion follows.

\begin{cor}[{Theorem 2.1 \cite{namioka1972}}]\label{c:ma:2}
Let $G$ be a metrizable CHART group.
Then $G$ is a topological group.
\end{cor}

From Theorem \ref{t:seqcomp:1} and \ref{t:ma:2}
the following assertion follows.

\begin{cor}\label{c:ma:3}
Suppose that $MA(\tau)$ is satisfied for each $\tau<\gs$.
Let $G$ be a sequentially compact CHART group.
Then $G$ is a topological group.
\end{cor}

The author thanks the referee for useful comments. 

\bibliographystyle{elsarticle-num}
\bibliography{mchart}

\begin{thebibliography}{10}
\expandafter\ifx\csname url\endcsname\relax
  \def\url#1{\texttt{#1}}\fi
\expandafter\ifx\csname urlprefix\endcsname\relax\def\urlprefix{URL }\fi
\expandafter\ifx\csname href\endcsname\relax
  \def\href#1#2{#2} \def\path#1{#1}\fi

\bibitem{Montgomery1936}
D.~{Montgomery}, {Continuity in topological groups}, {Bull. Am. Math. Soc.} 42
  (1936) 879--882.
\newblock \href {https://doi.org/10.1090/S0002-9904-1936-06456-6}
  {\path{doi:10.1090/S0002-9904-1936-06456-6}}.

\bibitem{Arens1946}
R.~Arens, \href{http://www.jstor.org/stable/2371787}{Topologies for
  homeomorphism groups}, American Journal of Mathematics 68~(4) (1946)
  593--610.
\newline\urlprefix\url{http://www.jstor.org/stable/2371787}

\bibitem{Ellis1957}
R.~{Ellis}, {A note on the continuity of the inverse}, {Proc. Am. Math. Soc.} 8
  (1957) 372--373.
\newblock \href {https://doi.org/10.2307/2033747} {\path{doi:10.2307/2033747}}.

\bibitem{Ellis1957-2}
R.~Ellis, Locally compact transformation groups, Duke Mathematical Journal
  24~(2) (1957) 119--125.

\bibitem{Ellis1960}
R.~Ellis, A semigroup associated with a transformation group, Transactions of
  the American Mathematical Society 94~(2) (1960) 272--281.

\bibitem{namioka1972}
I.~Namioka, Right topological groups, distal flows, and a fixed-point theorem,
  Mathematical systems theory 6~(1) (1972) 193--209.

\bibitem{MoorsNamioka2013}
W.~B. Moors, I.~Namioka, Furstenberg’s structure theorem via chart groups,
  Ergodic Theory and Dynamical Systems 33~(3) (2013) 954--968.

\bibitem{GlasnerMegrelishvili2013}
E.~Glasner, M.~Megrelishvili, Banach representations and affine
  compactifications of dynamical systems, in: Asymptotic geometric analysis,
  Springer, 2013, pp. 75--144.

\bibitem{rezn2022-1}
E.~Reznichenko, \href{https://arxiv.org/abs/2205.06316}{Continuity in right
  semitopological groups} (2022).
\newblock \href {https://doi.org/10.48550/ARXIV.2205.06316}
  {\path{doi:10.48550/ARXIV.2205.06316}}.
\newline\urlprefix\url{https://arxiv.org/abs/2205.06316}

\bibitem{ruppert1975}
W.~Ruppert, Uber kompakte rechtstopologische gruppen mit gleichgradig stetigen
  linkstranslationen, Sitz. ber. d. Osterr. Akad. d. Wiss. Math.-naturw. Kl.
  184 (1975) 159--169.

\bibitem{Milnes1979}
P.~Milnes, Continuity properties of compact right topological groups, in:
  Mathematical Proceedings of the Cambridge Philosophical Society, Vol.~86,
  Cambridge University Press, 1979, pp. 427--435.

\bibitem{Moors2016}
W.~B. Moors, Fragmentable mappings and chart groups, Fundamenta Mathematicae
  234 (2016) 191--200.

\bibitem{at2009}
A.~Arhangel'skii, M.~Tkachenko,
  \href{https://doi.org/10.2991/978-94-91216-35-0}{Topological Groups and
  Related Structures}, Atlantis Press, 2008.
\newblock \href {https://doi.org/10.2991/978-94-91216-35-0}
  {\path{doi:10.2991/978-94-91216-35-0}}.
\newline\urlprefix\url{https://doi.org/10.2991/978-94-91216-35-0}

\bibitem{Hodel1984handbook}
R.~Hodel, Cardinal functions i, in: Handbook of set-theoretic topology,
  Elsevier, 1984, pp. 1--61.

\bibitem{arh-rezn2005}
A.~Arhangel'skii, E.~Reznichenko, Paratopological and semitopological groups
  versus topological groups, Topology and its Applications 151~(1-3) (2005)
  107--119.

\bibitem{vanDouwen1984handbook}
E.~K. Van~Douwen, The integers and topology, in: Handbook of set-theoretic
  topology, Elsevier, 1984, pp. 111--167.

\bibitem{Vaughan1984handbook}
J.~E. Vaughan, Countably compact and sequentially compact spaces, in: Handbook
  of set-theoretic topology, Elsevier, 1984, pp. 569--602.

\bibitem{Shapirovskii1980}
B.~Shapirovskii, {Maps onto Tikhonov cubes}, Russian Mathematical Surveys
  35~(3) (1980) 145.

\bibitem{rezn2022gbd}
E.~{Reznichenko}, {Generalization of Baire spaces using diagonal}, arXiv
  e-prints (2022) 35\href {http://arxiv.org/abs/2203.09389}
  {\path{arXiv:2203.09389}}, \href
  {https://doi.org/https://doi.org/10.48550/arXiv.2203.09389}
  {\path{doi:https://doi.org/10.48550/arXiv.2203.09389}}.

\bibitem{MilnesPym1992}
P.~Milnes, J.~Pym,
  \href{https://www.scopus.com/inward/record.uri?eid=2-s2.0-84968480767&doi=10.1090%2fS0002-9939-1992-1065088-1&partnerID=40&md5=8eba4cf3993aee61498e0845ba31cf4d}{Haar
  measure for compact right topological groups}, Proceedings of the American
  Mathematical Society 114~(2) (1992) 387 – 393, cited by: 15; All Open
  Access, Bronze Open Access, Green Open Access.
\newblock \href {https://doi.org/10.1090/S0002-9939-1992-1065088-1}
  {\path{doi:10.1090/S0002-9939-1992-1065088-1}}.
\newline\urlprefix\url{https://www.scopus.com/inward/record.uri?eid=2-s2.0-84968480767&doi=10.1090%2fS0002-9939-1992-1065088-1&partnerID=40&md5=8eba4cf3993aee61498e0845ba31cf4d}

\bibitem{Moors2015}
W.~B. Moors, Invariant means on chart groups, Khayyam Journal of Mathematics
  1~(1) (2015) 36--44.

\end{thebibliography}
\end{document}